\documentclass[12pt]{article}

\usepackage{amsmath}   
\usepackage{amssymb}                  
\usepackage{amsthm}
\newcommand{\F}{{\mathbb F}}

\newcommand{\B}{{\mathcal B}}

\newcommand{\Mat}{{\mathrm{Mat}}}
\newcommand{\id}{{\mathrm{id}}}

\newcommand{\dt}{\textbf}

\newcommand{\be}{\begin{eqnarray*}}
\newcommand{\ee}{\end{eqnarray*}}

\newtheorem{defn}{Definition}[section]

\newtheorem{lem}{Lemma}[section]
\newtheorem{thm}{Theorem}[section]
\newtheorem{cor}{Corollary}[section]

\title{On finite congruence-simple semirings }
\author{  Chris Monico \\
          Department of Mathematics and Statistics\\ 
          Texas Tech University \\
          Lubbock, TX 79409-1042 \\
          \small{cmonico@math.ttu.edu}}
                                                                                                                                    
\begin{document}\maketitle
\thispagestyle{empty}

%%%%%%%%%%%%%%%%%%%%%%%%%%%%%%%%%%%%%%%%%%%%%%%%%%%%%%%%%%%%%%%%%% 
\begin{abstract}
  In this paper, we describe finite, additively commutative,
  congruence simple semirings. The main result is that the only such
  semirings are those of order 2, zero-multiplication rings of prime
  order, matrix rings over finite fields, ones with trivial addition
  and those that are additively idempotent.
\end{abstract}

%\begin{keyword} semirings \sep congruence simple \sep
%  congruence free
%\end{keyword}

%\end{frontmatter}
%%%%%%%%%%%%%%%%%%%%%%%%%%%%%%%%%%%%%%%%%%%%%%

\section{Introduction to semirings}
The notion of semirings seems to have first appeared
in the literature in a 1934 paper by Vandiver \cite{va34}.
Though the concept of a semiring might seem a bit strange
and unmotivated, additively commutative semirings arise 
naturally as the endomorphisms of commutative semigroups. 
Furthermore, every such semiring is
isomorphic to a sub-semiring of such endomorphisms \cite{he96a}.
For a more thorough introduction to semirings and a large collection
of references, the reader is referred to \cite{he96a}.

\begin{defn} A \dt{semiring} is a nonempty set $S$ together with
  two associative operations, $+$ and $\cdot$, such that for all
  $a,b,c\in S$, $a\cdot(b+c) = a\cdot b + a\cdot c$ and  
    $(a+b)\cdot c = a\cdot c + b\cdot c$.
\end{defn}
  A semiring is called \dt{additively [multiplicatively] commutative} 
  if $(S,+)$
  [$(S,\cdot)$] is commutative. If both $(S, +)$ and
  $(S,\cdot)$ are commutative, $S$ is simply called
  \dt{commutative}. 

The classification of finitely generated c-simple (see Definition
\ref{defsimple}) commutative semirings has only recently been
given in \cite{ba01}. In this paper, we progress toward a classification
of the class of finite c-simple semirings which are only additively commutative.
The main result will be given in Theorem \ref{thmMainch3}.

\begin{defn}
  An element $\alpha$ of a semiring is called 
  \dt{additively [multiplicatively] absorbing}
  if $\alpha + x = x+\alpha = \alpha$ [$\alpha \cdot x = x\cdot\alpha = \alpha$]
  for all $x\in S$. An element $\infty$ of a semiring
  is called an \dt{infinity} if it is both additively and
  multiplicatively absorbing.
\end{defn}

  Note that an additive identity in a semiring need not
  be multiplicatively absorbing. If, however, a semiring has a
  multiplicatively absorbing additive identity, we call it
  a \dt{zero}, and denote it by $0$.
  A semiring $S$ with additive identity $o$ is called 
  \dt{zero-sum free} if for all $a,b\in S$, 
  $a+b=o$ implies $a=b=o$.

\begin{defn}
  Let $S$ be a semiring and $\B\subseteq S$ a subset.
  Then $\B$ is called a \dt{bi-ideal} of $S$ if
  for all $b\in \B$ and $s\in S$,
  $ b+s, s+b, bs, sb \in \B.$
\end{defn}

\begin{defn} \label{defsimple}
  A \dt{congruence relation} on a semiring $S$ is
  an equivalence relation $\sim$ that also satisfies
\[
    x_1 \sim x_2  \Rightarrow
   \left\{
     \begin{array}{rcl}
        c+x_1 & \sim & c+x_2,   \\
        x_1+c & \sim & x_2+c,   \\
        cx_1 & \sim & cx_2,  \\
        x_1c & \sim & x_2c, 
      \end{array}
    \right. 
\]
  for all $x_1, x_2, c\in S$.
  A semiring $S$ that admits no congruence relations other than
  the trivial ones, $\id_S$ and $S\times S$, is said
  to be \dt{congruence-simple}, or \dt{c-simple}. 
\end{defn}

  Note that the trivial semiring of order 1 and every semiring of
  order 2 are congruence-simple. Also note that if $\B\subseteq S$
  is a bi-ideal then $\id_S\cup (\B\times \B)$ is a congruence relation.
  Thus, if $\B\subseteq S$ is a bi-ideal and $S$ is c-simple,
  then $\vert \B\vert = 1$ or $\B=S$.

The following theorem, due to Bashir, Hurt, Jan\v{c}a\v{r}\'{e}k, and Kepka
in \cite[Theorem 14.1]{ba01},
classifies finite c-simple commutative semirings.

\begin{thm} \label{thmBHJK}
  Let $S$ be a commutative, congruence-simple, finite semiring. 
  Then one of the following holds:
  \begin{enumerate}
    \item
      $S$ is isomorphic to one of the five semirings $T_1,\ldots, T_5$
       of order 2 defined in Table I.
    \item
      $S$ is a finite field.
    \item
      $S$ is a zero-multiplication ring of prime order.
    \item
      $S$ is isomorphic to $V(G)$ (defined below), for some
      finite abelian group $G$.
  \end{enumerate}
\end{thm}

  For a multiplicative abelian group $G$, set $V(G) = G\cup \{\infty \}$.
Extend the multiplication of $G$ to $V(G)$ by the rule 
$x\infty = \infty x = \infty$
for all $x\in V(G)$. Define an addition on $V(G)$ by the rules
$x+x = x$, $x+y = \infty$ for all $x,y\in V(G)$ with $x\ne y$. \\

%\noindent
\begin{center}
%\begin{table}[hbt] 
Table I\\
COMMUTATIVE SEMIRINGS OF ORDER TWO \\[12pt]
\begin{tabular}{|ll|ll|}
 \hline
  & & & \\
  \begin{tabular}{r|cc}
    $(T_1,+)$ & 0 & 1 \\
    \hline 
    0 & 0 & 0 \\
    1 & 0 & 0
  \end{tabular} \hspace{0.8 cm} &
  \begin{tabular}{r|cc}
    $\cdot$ & 0 & 1 \\
    \hline 
    0 & 0 & 0 \\
    1 & 0 & 0
  \end{tabular} &
  \begin{tabular}{r|cc} 
    $(T_2,+)$ & 0 & 1 \\
    \hline 
    0 & 0 & 0 \\
    1 & 0 & 0
  \end{tabular} \hspace{0.8 cm} &
  \begin{tabular}{r|cc}
    $\cdot$ & 0 & 1 \\
    \hline 
    0 & 0 & 0 \\
    1 & 0 & 1
  \end{tabular} \\
  & & & \\
 \hline
  & & & \\
  \begin{tabular}{r|cc}
    $(T_3,+)$ & 0 & 1 \\
    \hline 
    0 & 0 & 0 \\
    1 & 0 & 1
  \end{tabular} \hspace{0.8 cm} &
  \begin{tabular}{r|cc}
    $\cdot$ & 0 & 1 \\
    \hline 
    0 & 0 & 0 \\
    1 & 0 & 0
  \end{tabular} &
  \begin{tabular}{r|cc}
    $(T_4,+)$ & 0 & 1 \\
    \hline 
    0 & 0 & 0 \\
    1 & 0 & 1
  \end{tabular} \hspace{0.8 cm} &
  \begin{tabular}{r|cc}
    $\cdot$ & 0 & 1 \\
    \hline 
    0 & 1 & 1 \\
    1 & 1 & 1
  \end{tabular} \\
  & & & \\
 \hline
  & & & \\
  \begin{tabular}{r|cc}
    $(T_5,+)$ & 0 & 1 \\
    \hline 
    0 & 0 & 0 \\
    1 & 0 & 1
  \end{tabular} \hspace{0.8 cm} &
  \begin{tabular}{r|cc}
    $\cdot$ & 0 & 1 \\
    \hline 
    0 & 0 & 1 \\
    1 & 1 & 1
  \end{tabular} &
  \begin{tabular}{r|cc}
    $(T_6,+)$ & 0 & 1 \\
    \hline 
    0 & 0 & 0 \\
    1 & 0 & 1
  \end{tabular} \hspace{0.8 cm} &
  \begin{tabular}{r|cc}
    $\cdot$ & 0 & 1 \\
    \hline 
    0 & 0 & 0 \\
    1 & 0 & 1
  \end{tabular} \\
  & & & \\
 \hline
  & & & \\
  \begin{tabular}{r|cc}
    $(T_7,+)$ & 0 & 1 \\
    \hline 
    0 & 0 & 1 \\
    1 & 1 & 0
  \end{tabular} \hspace{0.8 cm} &
  \begin{tabular}{r|cc}
    $\cdot$ & 0 & 1 \\
    \hline 
    0 & 0 & 0 \\
    1 & 0 & 0
  \end{tabular} &
  \begin{tabular}{r|cc}
    $(T_8,+)$ & 0 & 1 \\
    \hline 
    0 & 0 & 1 \\
    1 & 1 & 0
  \end{tabular} \hspace{0.8 cm} &
  \begin{tabular}{r|cc}
    $\cdot$ & 0 & 1 \\
    \hline 
    0 & 0 & 0 \\
    1 & 0 & 1
  \end{tabular} \\
  & & & \\
 \hline
\end{tabular}
%\end{table}
\end{center}

We first note that a complete classification
up to isomorphism of finite, additively commutative, c-simple semirings is 
probably not possible.
To see this, note that $V(G)$ is c-simple for any finite
group $G$. Furthermore, if $G_1$ and $G_2$ are
two non-isomorphic groups, then $V(G_1)$ and $V(G_2)$ are
non-isomorphic semirings. Thus a classification of finite, additively
commutative, c-simple semirings up to isomorphism would require
a classification of finite groups.

%%%%%%%%%%%%%%%%%%%%%%%%%%%%%%%%%%%%%%%%%%%%%%%%%%%%%%%%%%%%%%%%%%%%

\section{Basic results}
The goal of this section is to derive some basic structure
information for finite, additively commutative, c-simple
semirings.

\begin{lem} \label{lemcen1}
  Let $S$ be a finite, additively commutative, c-simple semiring.
  If the multiplication table of $S$ has two
  identical rows [columns], then one of the following holds.
  \begin{enumerate}
    \item
      There exists $c\in S$ such that $xy=c$ for
      all $x,y\in S$.
    \item
      $\vert S\vert = 2$.
  \end{enumerate}
\end{lem}
\begin{proof}
  Observe that the relation $\sim$ defined by
  \[
    x\sim y \hspace{0.5cm}\mbox{ if }\hspace{0.5cm} xz = yz 
    \;\mbox{ for all } z\in S
  \]
  is a congruence
  relation. By assumption, there exist $r_1\ne r_2$
  such that $r_1z = r_2z$ for all $z\in S$ so
  $\sim\,= S\times S$. Thus 
  \begin{equation} \label{eqRMult}
    xz=yz \hspace{0.5cm}\mbox{ for all } x,y,z\in S.
  \end{equation}

  Suppose that $(S,\cdot)$ is not left-cancellative. Then
  there exist $a,b,c,d\in S$ such that $da=db=c$ and
  $a\ne b$. But $xa = ya, xb=yb$ for all $x,y\in S$.
  Hence $da = ya, db=yb$ and so $ya=yb=c$ for all 
  $y\in S$.
  Consider now the congruence relation $\approx$
  defined by 
  \[
    x\approx y\hspace{0.5cm}\mbox{ if } \hspace{0.5cm} zx=zy 
    \;\mbox{ for all } z\in S.
  \]
  Since $a\ne b$ and $a\approx b$, it follows that
  $\approx \,= S\times S$, whence $zx=zy$ for all $x,y,z\in S$.
  Then for all $x,y\in S$ we have
  $xy = xa = da = c$.

  Suppose now that $(S,\cdot)$ is left-cancellative.
  Fix $x\in S$ and let $z=x^2$. Then $xz=zx$.
  But $yz = xz$ and $yx = zx$ for all
  $y\in S$, so $yz = yx$. By 
  left-cancellation, $x^2 = z = x$, so $S$ is multiplicatively
  idempotent.  Furthermore, for all $w\in S$
  $w+w = w^2+w^2 = (w+w)w = w^2 = w$, so $S$ is additively
  idempotent. We will now show, by contradiction, that
  $\vert S\vert \le 2$.

  Suppose $\vert S\vert = n>2$. 
  For each nonempty subset $A\subseteq S$ let
  \[ 
    \sigma_A = \sum_{x\in A}x
  \]
  and $\sigma = \sigma_S$.
  Suppose that $A\subset S$ with $\vert A\vert = n-1$. Consider the
  relation 
  $\sim \,= \id_S\cup \{ (\sigma_A, \sigma),
                                (\sigma, \sigma_A)\}$.
  Clearly $\sim$ is an equivalence relation. Since $(S,\cdot)$
  is idempotent, Equation \ref{eqRMult} implies that
  for each $c\in S$
  \[
    c\sigma_A = \sigma_A\sigma_A = \sigma_A
    \hspace{0.5cm}\mbox{ and }\hspace{0.5cm}
    c\sigma = \sigma\sigma = \sigma.
  \]
  Thus, $c\sigma_A \sim c\sigma$. Similarly,
  \[
    \sigma_A c = c^2 = c
    \hspace{0.5cm}\mbox{ and }\hspace{0.5cm}
    \sigma c = c^2 = c
  \]
  so that $\sigma_A c \sim \sigma c$.
  Since $(S,+)$ is idempotent, $\sigma+c=\sigma$ and
  \[
    \sigma_A + c = \left\{
      \begin{array}{ll}
        \sigma_A, & \mbox{ if $c\in A$}, \\
        \sigma,   & \mbox{ otherwise }.
      \end{array} \right.
  \]
  Thus $\sim$ is a congruence relation. Since 
  $\vert S\vert > 2$, it must be the case that
  $\sim \,= \id_S$, so $\sigma_A = \sigma$ for all
  proper $A\subset S$ with $\vert A\vert = n-1$.
  
  By induction, we will now show that $\sigma_A = \sigma$
  for any nonempty subset $A\subseteq S$. Suppose this
  is known to hold for all $A$ with 
  $\vert A\vert = k \ge 2$. Let $A\subset S$ with
  $\vert A\vert = k-1$ and again consider the relation
  \[
    \sim \,= \id_S\cup 
        \{ (\sigma_A, \sigma), (\sigma, \sigma_A)\}.
  \]
  As above, $\sim$ is a multiplicative equivalence relation.
  Furthermore
  \[
    \sigma_A + c = \left\{
      \begin{array}{ll}
        \sigma_A, & \mbox{ if $c\in A$}, \\
        \sigma_{A\cup\{c\}},   & \mbox{ otherwise }.
      \end{array} \right.
  \]
  But $c\not\in A$ implies $\vert A\cup\{c\}\vert = k$,
  so $\sigma_{A\cup\{c\}} = \sigma$ by the inductive assumption.
  Thus $\sim$ is again a congruence relation. Since
  $\sim \,\ne S\times S$, it follows that $\sim \,= \id_S$, so
  $\sigma_A = \sigma$.

  In particular, this shows that for each $w\in S$,
  $ w = \sigma_{\{w\}} = \sigma$,
  a contradiction. Thus $\vert S\vert = 2$.

  It only remains to see that the same statement holds if
  ``rows'' is replaced by ``columns''.
  If $S$ has two identical columns, consider the
  reciprocal semiring $(S', +, \otimes)$ defined by
  $(S', +) = (S, +)$ and $x\otimes y = yx$.
  This semiring is c-simple and has two identical rows
  so the above argument applies. 
\end{proof}

\begin{lem} \label{lemcen4}
  Let $S$ be a finite, additively commutative, c-simple
  semiring. Then one of the following holds.
  \begin{itemize}
    \item
      $(S,+)$ is a group, hence $(S,+,\cdot)$ is a ring.
    \item
      $S$ has an additively absorbing element $\alpha$.
  \end{itemize}
\end{lem}
\begin{proof}
  Consider the relation $\sim$ defined by
  \[
    x\sim y \hspace{0.5cm} \mbox{ if }\hspace{0.5cm}
     x+t=y+t \,\mbox{ for some } t\in S.
  \]
  It is easy to see that $\sim$ is a congruence relation.
  If $\sim \,= \id_S$, then $(S,+)$ is cancellative, hence
  a group. It follows easily that $(S,+,\cdot)$ is a ring.
  On the other hand, suppose $\sim \,= S\times S$.
  Then for all $x,y\in S$ there exists $t_{x,y}\in S$
  such that $x+t_{x,y} = y + t_{x,y}$. Set
  \[
    \sigma = \sum_{x\in S}x 
      \hspace{0.75cm}\mbox{ and }\hspace{0.75cm}
    \alpha = \sigma+\sigma.
  \]
  For $x,y\in S$ there exists $\sigma'\in S$ such that
  $\sigma = t_{x,y} + \sigma'$. Then
  \[
    x+\sigma = x+t_{x,y}+\sigma' = y + t_{x,y} + \sigma' = y+\sigma.
  \]
  In particular, $x+\sigma = \sigma+\sigma$ for all $x\in S$.
  Thus, for all $x\in S$
  \[
    x+\alpha = x+\sigma+\sigma = (\sigma+\sigma)+\sigma
             = \sigma+\sigma = \alpha.
  \]

\end{proof}

\begin{thm} \label{thmssr1}
  Let $S$ be a finite, additively commutative, c-simple
  semiring. Then one of the following holds.
  \begin{itemize}
    \item
      $(S,+, \cdot)$ is a ring.
    \item
      $S$ has an infinity.
    \item
      $S$ is additively idempotent.
  \end{itemize}
\end{thm}
\begin{proof}
  With respect to Lemma \ref{lemcen4}, one may assume that there
  is an additively absorbing element $\alpha\in S$.
  Consider the relation $T$ defined by
  \[
    xTy 
      \hspace{0.5cm}\mbox{ if }\hspace{0.5cm} 2x = 2y.
  \]
  Then $T$ is a congruence relation, whence $T=\id_S$ or
  $T=S\times S$.

  \noindent
  {\em Case I:} Suppose $T=S\times S$.

    Then for all $x\in S$,
  $x+x = \alpha + \alpha = \alpha$. Thus,
  $x\alpha = x(\alpha + \alpha) = x\alpha + x\alpha = \alpha$.
  Similarly, $\alpha x = \alpha$ so $\alpha$ is an infinity.

  \noindent
  {\em Case II:} Suppose $T=\id_S$.

  Consider the congruence relation $\sim$
  defined by $x\sim y$ if there exist $u,v\in S\cup\{o\}$ and
  $i\ge 0$ such that
  \be
    2^ix &=& y+u, \\
    2^iy &=& x+v.
  \ee
  Then $2(2x) = (x)+3x$ and $2(x) = (2x) + o$, so $x\sim 2x$
  for all $x\in S$. If $\sim \,=\id_S$, then $x=2x$ for all
  $x\in S$, whence $(S,+)$ is idempotent.
  Suppose now that $\sim \,= S\times S$ and let $x\in S$.
  Then $x\alpha\sim \alpha$, so there exists $v\in S\cup\{o\}$ and
  $i\ge 0$ such that $2^ix\alpha = \alpha+v = \alpha$.
  Then
  \[
    x\alpha = x(2^i\alpha) = 2^ix\alpha = \alpha,
  \]
  so $x\alpha = \alpha$. Similarly, $\alpha x = \alpha$
  so $\alpha$ is an infinity.
\end{proof}

\begin{cor} \label{cor0cent} 
  If $S$ is a finite, additively commutative, 
  c-simple semiring with zero then one of the following holds.
  \begin{itemize} 
    \item 
      $S\cong \Mat_n(\F_q)$ for some $n\ge 1$ and some 
      finite field $\F_q$.
    \item 
      $S$ is a zero-multiplication ring ($S^2=\{0\}$) 
      of prime order. 
    \item 
      $S$ is additively idempotent. 
  \end{itemize} 
\end{cor}

%%%%%%%%%%%%%%%%%%%%%%%%%%%%%%%%%%%%%%%%%%%%%%%%%%%%%%%%%%%%%%%%% 
\section{The $\infty$ case} 
In this section, we show that a finite, additively commutative, 
c-simple semiring with $\infty$ is either additively 
idempotent, has trivial addition, or has order 2. 
 
\begin{lem} \label{leminf1} 
  Let $S$ be a finite, additively commutative, c-simple semiring 
  with $\infty$ and $\vert S\vert > 2$. Then one of the following
  holds
  \begin{enumerate}
    \item
      $S$ is additively idempotent.
    \item
      $S+S=\{\infty\}$ and $(S,\cdot)$ is a congruence-free semigroup.
  \end{enumerate}
\end{lem} 
\begin{proof} 
  Consider the congruence relation defined by
  \[ 
    xTy \hspace{0.5cm}\mbox{ if }\hspace{0.5cm} 
      2x=2y. 
  \] 
  {\em Case I:} \hspace{0.25cm}$T=\id_S$. 

    Then $2x=2y$ iff $x=y$. Set $x\sim y$ if there exists $i\ge 0$ and  
    $u,v\in S\cup\{o\}$ such that 
  \be 
    2^ix &=& y+u, \\ 
    2^iy &=& x+v. 
  \ee 
  Then $\sim$ is a congruence relation and $x\sim 2x$ 
  for all $x\in S$. But $x\not\sim \infty$ for $x\ne\infty$, 
  so $\sim \,\ne S\times S$. Thus, $\sim \,=\id_S$, and so 
  $S$ is additively idempotent.\\ 
  \noindent 
  {\em Case II:}\hspace{0.25cm} $T=S\times S$.

  Then $x+x=\infty$ for all $x\in S$. For 
  $\emptyset\ne A\subseteq S$, let  
  \[ 
    \sigma_A = \sum_{x\in A}x .
  \] 
  Let $N=\vert S\vert$ 
  and suppose that $\vert A\vert = N-1$. Then for every $c\in S$,
  $\sigma_A + c =\infty$, since $c\in A$, $c=\infty$, or $\sigma_A=\infty$.
  Furthermore,
  \[ 
    c\sigma_A =\sum_{x\in A}cx = 
    \left\{ \begin{array}{ll} 
      \infty, & \mbox{ if $cx_1=cx_2$ for some distinct $x_1,x_2\in A$}, \\ 
      \sigma_A, & \mbox{ otherwise}.
    \end{array}\right. 
  \] 
  Similarly, $\sigma_A c =\infty$ or $\sigma_A c=\sigma_A$. 
  Thus, $\B =\{\sigma_A \,\vert\, A\subset S \mbox{ with }  
               \vert A\vert = N-1\}$ is a bi-ideal.
  Furthermore, $\infty\in A$ implies $\sigma_A=\infty$. 
  Thus, $\vert\B\vert\le 2$ and  so $\B=S\Rightarrow \vert S\vert = 2$,
  a contradiction.
  Thus $\B = \{\infty\}$, so $\sigma_A=\infty$ for all 
  $A\subset S$ with $\vert A\vert = N-1$. 
 
  By induction, we will show that $\sigma_A=\infty$ for all 
  $A\subset S$ with $\vert A\vert = 2$. 
  Assume $\sigma_A=\infty$ for all $A\subset S$ with 
  $\vert A\vert = k+1 > 2$. 
 
  Suppose now that $A\subset S$ with $\vert A\vert = k\ge 2$. 
  Then for $c\in S$, 
  \[ 
    \sigma_A + c = \left\{ \begin{array}{ll} 
      \infty, & \mbox{ if } c\in A, \\ 
      \sigma_{A\cup\{c\}}, & \mbox{ otherwise }. 
    \end{array}\right. 
  \] 
  By assumption, if $c\not\in A$ then $\sigma_{A\cup\{c\}} = \infty$, 
  so $\sigma_A + c = \infty$ for all $c\in S$. 
  Also
  \[ 
    c\sigma_A =\sum_{x\in A}cx = 
    \left\{ \begin{array}{ll} 
      \infty, & \mbox{ if $cx_1=cx_2$ for some distinct $x_1,x_2\in A$}, \\ 
      \sigma_B, & \mbox{ for some $\vert B\vert=k$ otherwise}. 
    \end{array}\right. 
  \] 
  The same is easily seen to hold for $\sigma_Ac$. Observe that 
  $\sigma_X=\infty$ for some $X\subset S$ with $\vert X\vert = k$, 
  so 
  \[ 
    \B = \{\sigma_A \,\vert\, A\subset S \mbox{ with } \vert A\vert = k\} 
  \] 
  is a bi-ideal of $S$.\\ 
  \noindent 
  {\em Case (i):} \hspace{0.25cm}$\B=\{\infty\}$.

  Then $\sigma_A=\infty$ for all $A\subset S$ with $\vert A\vert = k$, 
  so we may apply the induction and conclude that  
  $\sigma_A =\infty$ for all $A\subset S$ with $\vert A\vert = 2$. 
  Thus, $x+y=\infty$ for all $x,y\in S$.\\ 
  \noindent 
  {\em Case (ii):} \hspace{0.25cm}$\B=S$.  

  We will show directly that $x+y=\infty$ for all $x,y\in S$. 
  By assumption this holds for $x=y$, so suppose $x\ne y$. 
  Then there exist $A_1,A_2\subset S$ with  
  $\vert A_1\vert = \vert A_2\vert = k$ and  
  $\sigma_{A_1} = x$, $\sigma_{A_2}=y$. 
  \be 
    A_1\cap A_2\ne\emptyset & \Rightarrow & x+y =  
             \sigma_{A_1} + \sigma_{A_2} =\infty. \\ 
    A_1\cap A_2=\emptyset & \Rightarrow & x+y =  
             \sigma_{A_1} + \sigma_{A_2} =\sigma_{A_1\cup A_2}.  
  \ee 
  But $\vert A_1\cup A_2\vert > k$. In particular, either 
  $\vert A_1\cup A_2\vert = k+1$ or there exist 
  $\emptyset\ne B_1, B_2\subset S$ with $\vert B_1\vert = k+1$, 
  $B_1\cap B_2 =\emptyset$ and $B_1\cup B_2 = A_1\cup A_2$. 
  By assumption, $\sigma_{B_1}=\infty$ and we have 
  \[ 
    x+y =\sigma_{A_1\cup A_2} = \sigma_{B_1\cup B_2}  
                              = \sigma_{B_1} + \sigma_{B_2} 
                              = \infty + \sigma_{B_2} = \infty.
  \] 
  Thus $x+y=\infty$ for all $x,y\in S$. Finally, note that
  since $S+S=\{\infty\}$, any nontrivial 
  congruence relation on $(S,\cdot)$
  is also a nontrivial congruence relation on $(S,+,\cdot)$, whence
  $(S,\cdot)$ is a congruence-free semigroup.
\end{proof} 

The following is Theorem 3.7.1 from \cite{ho95}.

\begin{thm} \label{thmHO}
  Let $I=\{1,2,\ldots, m\}$, $\Lambda=\{1,2,\ldots, n\}$, and
  $P=(p_{ij})$ be an $n\times m$ matrix of 1's and 0's such that no row or column is
  identically zero, no two rows are identical, and no two columns
  are identical. Let $S=(I\times\Lambda)\cup\{\infty \}$ and define
  a binary relation on $S$ by
  \[
    (i,\lambda)\cdot (j, \mu) =
    \left\{ \begin{array}{ll}
              (i,\mu) & \mbox{ if $p_{\lambda j} = 1$} \\
              \infty       & \mbox{ otherwise,}
            \end{array} \right.
  \]
  \[
    (i,\lambda)\cdot \infty = \infty\cdot (i,\lambda) = \infty\cdot \infty = \infty.
  \]
  Then $S$ is a congruence-free semigroup of order $mn+1$.
  Conversely, every finite congruence-free semigroup with
  an absorbing element is isomorphic to one of this kind.
\end{thm}

%%%%%%%%%%%%%%%%%%%%%%%%%%%%%%%%%%%%%%%%%%%%%%%%%%%%%%%%%%%%%%%%%%%% 
\section{Main theorem}

  \begin{thm} \label{thmMainch3}
    Let $S$ be a finite, additively commutative, congruence-simple semiring.
    Then one of the following holds:
    \begin{enumerate}
      \item
        $\vert S\vert \le 2$.
      \item
        $S\cong\Mat_n(\F_q)$ for some finite field $\F_q$
        and some $n\ge 1$.
      \item
        $S$ is a zero multiplication ring of prime order.
      \item
        $S$ is additively idempotent. 
      \item
        $(S,\cdot)$ is a semigroup as in Theorem \ref{thmHO} 
        with absorbing element $\infty\in S$ and $S+S=\{\infty\}$.
    \end{enumerate}
  \end{thm}
  \begin{proof}
    Apply Theorems \ref{thmssr1} and \ref{thmHO}, Lemma \ref{leminf1} and 
    Corollary \ref{cor0cent}. Also notice that if $(S,\cdot)$ is a semigroup
    as in Theorem \ref{thmHO}, and we define $S+S=\{\infty\}$, then
    $(S,+,\cdot)$ is necessarily congruence-free.
  \end{proof}
  Observe the similarity between this theorem and Theorem \ref{thmBHJK}.
  Recall that for a finite group $G$, $V(G)$ is a finite, additively
  commutative, c-simple semiring and is additively idempotent. So
  the semirings $V(G)$ do fall into the fourth case of Theorem
  \ref{thmMainch3}. Note also that for $n>1$, the matrix
  semiring $\Mat_n(V(G))$ is not c-simple. To see this, consider a
  matrix with all but one entry equal to infinity, and apply 
  Lemma \ref{lemcen1}.
  In view of this, it might be tempting to conjecture that the additively idempotent
  semirings are precisely those of the form $V(G)$. However, the
  semiring in Table III provides a counter-example to that conjecture.

\begin{table}[hbt] \label{tabssr3}
\begin{center}
Table III\\
A C-SIMPLE SEMIRING OF ORDER 3\\[12pt]
  \begin{tabular}{cc}
    \begin{tabular}{r|ccc}
      $+$ & a & 1 & b\\
      \hline
      a & a & 1 & b \\
      1 & 1 & 1 & b \\
      b & b & b & b
    \end{tabular} \hspace{0.8 cm} &
    \begin{tabular}{r|ccc}
      $\cdot$ & a & 1 & b\\
      \hline
      a & a & a & b \\
      1 & a & 1 & b \\
      b & a & b & b
    \end{tabular}
  \end{tabular}
 \end{center}
\end{table}
This semiring is additively idempotent yet has order 3
and is not of the form $V(G)$. At present, we have no strongly
supported conjecture for a meaningful description of the
semirings in the fourth case of Theorem \ref{thmMainch3},
though we do believe that some good description might be
possible.

%%%%%%%%%%%%%%%%%%%%%%%%%%%%%%%%%%%%%%%%%%%%%%%%%%%%%%%%%%%%%%%%%%%% 
\section{Acknowledgments}
  The author would like to thank the anonymous referee
  for his/her careful reading and valuable input which
  greatly increased the quality of this work.
  This research was supported
  by a fellowship from the Center for Applied Mathematics
  at the University of Notre Dame, and in part by NSF grant
  DMS-00-72383.

%\bibliographystyle{plain}
%\bibliography{huge}

\begin{thebibliography}{2}

\bibitem{ba01}
R.~El Bashir, J.~Hurt, A.~Jan\v{c}a\v{r}\'{e}k, and T.~Kepka.
\newblock Simple commutative semirings.
\newblock {\em Journal of Algebra}, 236 (2001), 277--306, doi:10.1006/jabr.2000.8483

\bibitem{he96a}
Udo Hebisch and Hanns~Joachim Weinert.
\newblock Semirings and Semifields.
\newblock in ``Handbook of Algebra'', Vol. 1., Elsevier Science B.V., Amsterdam, 1996.

\bibitem{ho95}
J. Howie.
\newblock Fundamentals of Semigroup Theory.
\newblock Oxford University Press, NY, 1995.

\bibitem{va34}
H.S. Vandiver.
\newblock Note on a simple type of algebra in which the cancellation law of
  addition does not hold.
\newblock {\em Bulletin of the American Mathematical Society}, 40 (1934), 916--920.

\end{thebibliography}

\end{document}